\documentclass[12pt]{amsart}
\usepackage{amsmath}
\usepackage{amsfonts}
\usepackage{amssymb}
\usepackage{amsthm}
\usepackage{amstext,amscd}

\theoremstyle{plain}
\newtheorem{theorem}{Theorem}[section]
\newtheorem*{acknowledgements}{Acknowledgements}

\newtheorem{lemma}[theorem]{Lemma}

\newtheorem{proposition}[theorem]{Proposition}
\newtheorem{remark}{Remark}

\numberwithin{equation}{section}

\begin{document}
\title{On the Generalized Volume Conjecture and Regulator}
\author{Weiping Li}
\address{Department of Mathematics\\
         Oklahoma State University\\
         Stillwater, OK 74078}
\email{wli@math.okstate.edu}
\author{Qingxue Wang}
\address{Chern Institute of Mathematics\\
         Nankai University\\
         Tianjin 300071\\
         P.R. China}
\email{qingxue\_wang@yahoo.com.cn}

\begin{abstract}
In this paper, by using the regulator map of Beilinson-Deligne on a
curve, we show that the quantization condition posed by Gukov is
true for the $SL_2(\mathbb{C})$ character variety of the hyperbolic
knot in $S^3$. Furthermore, we prove that the corresponding
$\mathbb{C}^{*}$-valued closed $1$-form is a secondary
characteristic class (Chern-Simons) arising from the vanishing first
Chern class of the flat line bundle over the smooth part of the
character variety, where the flat line bundle is the pullback of the
universal Heisenberg line bundle over $\mathbb{C}^{*}\times
\mathbb{C}^{*}$. Based on this result, we give a reformulation of
Gukov's generalized volume conjecture from a motivic perspective.
\end{abstract}

\maketitle

\section{Introduction}
It is a very important question in knot theory to find the geometric
and topological interpretation of the Jones polynomial of a knot.
Observed by Kashaev \cite{Kas}, H. Murakami and J. Murakami
\cite{MM}, the asymptotic rate of (N-colored) Jones polynomial is
related to the volume of the hyperbolic knot complement. It is known
as the Volume conjecture. Following Witten's $SU(2)$ topological
quantum field theory, Gukov \cite{Guk} proposed a complex version of
Chern-Simons theory and generalized the volume conjecture to a
$\mathbb{C}^{*}$- parametrized version with parameter lying on the
zero locus of the $A$-polynomial of the knot in $S^3$.

In this paper, we prove that the quantization condition posed by
Gukov \cite[Page~597]{Guk} is true for hyperbolic knots in $S^3$.
The key ingredient of the proof is the construction of the regulator
map of an algebraic curve studied by Beilinson, Bloch, Deligne and
many others (\cite{Bei,Bl,De,Ram}). Let $K$ be a hyperbolic knot in
$S^3$. For each irreducible component $Y$ of the zero locus of the
$A$-polynimial $A(l,m)$ of $K$, we show that the symbol $\{l,m\}\in
K_2(\mathbb{C}(Y))$ is a torsion. Associated to $\{l,m\}$, there is
a cohomology class $r(l,m)$ in $H^1(Y_h, \mathbb{C}^{*})$, where
$Y_h$ is some open Riemann surface. The detail is given in Section
$3$. As Deligne noted, $H^1(Y_h, \mathbb{C}^{*})$ is the group of
isomorphism classes of flat line bundles over $Y_h$. Thus, our class
$r(l,m)$ corresponds to a flat line bundle. Moreover, the line
bundle $r(l,m)$ can be constructed explicitly as the pullback of the
universal Heisenberg line bundle over $\mathbb{C}^{*}\times
\mathbb{C}^{*}$, see \cite{Bl,Ram}. We then derive a closed $1$-form
from it and show that this \emph{$1$-form} is the Chern-Simons class
of the \emph{first} Chern class $C_1$ of $r(l,m)$. Note that our
Chern-Simons class is not the usual Chern-Simons class as a closed
\emph{$3$-form} of the \emph{second} Chern class for a
$3$-dimensional manifold. We also reformulate the generalized volume
conjecture via this closed $1$-form Chern-Simons class. Our proof is
motivic in nature.

The paper is organized as follows. In section $2$, we introduce the
notations used in the paper. In section $3$, we discuss the
generalized volume conjecture and the regulator of a curve, then we
prove the theorem about the quantization condition and give a
reformulation of the generalized volume conjecture. In the end, we
remark the motivic aspect of the proof.

\begin{acknowledgements} We would like to thank Professor Alexander
Goncharov and Professor Xiao-Song Lin for their helpful comments and
suggestions on an earlier version of the paper.
\end{acknowledgements}

\section{Terminology and Notation}
\subsection{}
Let $K$ be a knot in $S^3$ and $M_K$ its complement. That is,
$M_K=S^3-N_K$ where $N_K$ is the open tubular neighborhood of $K$ in
$S^3$. $M_K$ is a compact $3$-manifold with boundary $\partial
M_K=T^2$ a torus. Denote by
$R(M_K)=\text{Hom}(\pi_1(M_K),SL_2(\mathbb{C}))$ and $R(\partial
M_K)=\text{Hom}(\pi_1(\partial M_K),SL_2(\mathbb{C}))$. It is known
that they are affine algebraic sets over the complex numbers
$\mathbb{C}$ and so are the corresponding character varieties
$X(M_K)$ and $X(\partial M_K)$ (See \cite{CS}). We also have the
canonical surjective morphisms $t:R(M_K)\longrightarrow X(M_K)$ and
$t:R(\partial M_K)\longrightarrow X(\partial M_K)$ which map a
representation to its character. The natural homomorphism $i:
\pi_1(\partial M_K)\longrightarrow \pi_1(M_K)$ induces the
restriction maps $r: X(M_K)\longrightarrow X(\partial M_K)$ and $r:
R(M_K)\longrightarrow R(\partial M_K)$.

\subsection{}
Since $\pi_1(\partial M_K)=\mathbb{Z}\oplus \mathbb{Z}$, we shall
fix two oriented simple curves $\mu$ and $\lambda$ as its
generators. They are called the meridian and longitude respectively.
Let $R_D$ be the subvariety of $R(\partial M_K)$ consisting of the
diagonal representations. Then $R_D$ is isomorphic to
$\mathbb{C}^{*} \times \mathbb{C}^{*}$. Indeed, for $\rho \in R_D$,
we obtain
\begin{equation*}
\rho(\lambda)=
  \left[
  \begin{matrix}
    l & 0\\
    0 & l^{-1}
  \end{matrix}
  \right] \, \text{and} \, \,
\rho(\mu)=
  \left[
  \begin{matrix}
    m & 0\\
    0 & m^{-1}
  \end{matrix}
  \right],
\end{equation*}
then we assign the pair $(l,m)$ to $\rho$. Clearly this is an
isomorphism. We shall denote by $t_{D}$ the restriction of the
morphism $t:R(\partial M_K)\longrightarrow X(\partial M_K)$ on
$R_D$.

\subsection{}
Next we recall the definition of the $A$-polynomial of $K$ which was
introduced in \cite{CCGLS}. Denote by $X^{\prime}(M_K)$ the union of
the irreducible components $Y^{\prime}$ of $X(M_K)$ such that the
closure $\overline{r(Y^{\prime})}$ in $X(\partial M_K)$ is
$1$-dimensional. For each component $Z^{\prime}$ of
$X^{\prime}(M_K)$, denote by $Z$ the curve
$t_{D}^{-1}(\overline{r(Y^{\prime})})\subset R_D$.  We define $D_K$
to be the union of the curves $Z$ as $Z^{\prime}$ varies over all
components of $X^{\prime}(M_K)$. Via the above identification of
$R_D$ with $\mathbb{C}^{*} \times \mathbb{C}^{*}$, $D_K$ is a curve
in $\mathbb{C}^{*} \times \mathbb{C}^{*}$. Now by definition the
\emph{$A$-polynomial $A(l,m)$} of $K$ is the defining polynomial of
the closure of $D_K$ in $\mathbb{C}\times\mathbb{C}$.

From now on, we shall assume that $K$ is a hyperbolic knot. Denote
by $\rho_0: \pi_1(M_K)\longrightarrow PSL_2(\mathbb{C})$ the
discrete, faithful representation corresponding to the hyperbolic
structure on $M_K$. Note that $\rho_0$ can be lifted to a
$SL_2(\mathbb{C})$ representation. Moreover, there are exactly
$|H^{1}(M_K; \mathbb{Z}_2)=\mathbb{Z}_2|=2$ such lifts.

\section{$A$-polynomial, regulator and $K_2$ of a curve}
In this section, we briefly recall Gukov's formulation of the
generalized volume conjecture. Using the regulator map of a curve,
we show that the form $r(l,m)=\xi(l,m)+i \eta(l,m)$ over the
$1$-dimensional character variety $Y_h$ has exact imaginary part and
rational real part. This provides an affirmative answer to Gukov's
quantization over $Y_h$. Moreover, $dr(l,m)=\frac{dl}{l}\wedge
\frac{dm}{m}=2 \pi i C_1(L)=0$ justifies that $r(l,m)$ is the
Chern-Simons class from the first Chern class $C_1$.

\subsection{}
Let $\overline{D_K}$ be the zero locus of the $A$-polynomial
$A(l,m)$ in $\mathbb{C}^2$. Let $y_0\in D_K$ correspond to the
character of the representation of the hyperbolic structure on $M_K$
and $m(y_0)=1$. For a path $c$ in $\overline{D_K}$ with the initial
point $y_0$ and endpoint $(l,m)$, the following quantities are
defined in \cite[~(5.2)]{Guk}:
\begin{equation}\label{vol}
   Vol(l,m)=Vol(K)+2\int_{c}[-\log{|l|}\,d(\arg{m})+\log{|m|}\,d(\arg{l})],
\end{equation}
\begin{equation}\label{cs}
   CS(l,m)=CS(K)-\frac{1}{\pi^2} \int_{c}[\log{|m|}\,d
\log{|l|}+(\arg{l})\,d(\arg{m})],
\end{equation}
where $Vol(K)$ and $CS(K)$ are the volume and the Chern-Simons
invariant of the complete hyperbolic metric on $M_K$.

In \cite[~(5.12)]{Guk}, Gukov proposes his \emph{Generalized Volume
Conjecture}: for a fixed number $a$ and
$m=-\exp({i \pi a})$,\\
\begin{equation}\label{conj}
\lim_{N, k\rightarrow \infty; \frac{N}{k}=a}
\frac{\log{J_{N}(K,e^{2\pi i/k})}}{k}=\frac{1}{2 \pi}(Vol(l,m)+
i2\pi^2 CS(l,m)),
\end{equation}
where $J_{N}(K,q)$ is the $N$-colored Jones polynomial of $K$,
$Vol(l,m)$ and $CS(l,m)$ as in (\ref{vol}) and (\ref{cs}), are the
functions on the zero locus of the A-polynomial of the hyperbolic
knot $K$. Note that for $m=1$ and $a=1$, we get the usual
\emph{Volume Conjecture}:
\begin{equation}
\lim_{N\rightarrow \infty} \frac{\log{|J_{N}(K,e^{2\pi
i/N})|}}{N}=\frac{1}{2 \pi}Vol(K).
\end{equation}

Gukov's generalized volume conjecture links the Jones invariants of
knots with the topological and geometric invariants arising from the
character variety of the knot complement. It has received lot of
attention. But other than the verification of a few examples, there
is no essential mathematical evidence to support this interesting
conjecture. And there appeared some confusions in the study of this
conjecture (see [Mu], [Gu-Mu]). Understanding those terms in Gukov's
generalized volume conjecture (\ref{conj}) would be the first step.

\begin{remark}
In \cite{Mu}, Murakami reformulated a parametrized volume conjecture
where $2\pi i$ is replaced by $2\pi i +u$. He conjectured that
\[
 H(K,u)=(2\pi i +u)\lim_{N\rightarrow
 \infty}\frac{log(J_N(K,e^{\frac{2\pi i +u}{N}}))}{N}
\]
is analytic on some open subset of $\mathbb{C}$ and the volume
function $V(K,u)$ is given by
\[
  V(K,u)=\text{Im}(H(K,u))-\text{Re}(u)\cdot
  \text{Im}(\frac{dH(K,u)}{du}).
\]
This parametrization is different from Gukov's original one
($m=\exp(u)$, $l=\exp(2\frac{dH}{du}-2\pi i)$). It comes from a
different choice of polarization when $\text{Re}(u)\ne 0$.
Throughout this paper, we shall keep the same polarization as
Gukov's.
\end{remark}
For $Vol(l,m)$ in (\ref{vol}), it is understood that it measures the
change of volumes of the representations on the path $c$ in
$\overline{D_K}$, the zero locus of $A$-polynomial of the hyperbolic
knot. See \cite[Sect.~4.5]{CCGLS} and
\cite[Sect.~2]{Dun} for more detail.\\

For $CS(l,m)$ in (\ref{cs}), to our knowledge, it has not been
understood mathematically. It was derived from the point of view of
physics, see \cite[Sect.~3]{Guk}. On the other hand, since $M_K$ has
boundary a torus $T$, by \cite{RSW} and \cite{KK}, its $3$-form
Chern-Simons functional is only well-defined as a section of a
circle bundle over the gauge equivalence classes of $T$. By
\cite[Theorem~3.2, 2.7]{KK}, if $\chi_t$, $t\in [0,1]$ is a path of
characters of $SL_2(\mathbb{C})$ representations of $M_K$ and $z(t)$
is the Chern-Simons invariant of $\chi_t$, then:
\begin{equation}\label{kkf}
  \begin{aligned}
     z(1)z(0)^{-1}&=\text{exp}(2\pi i(\int_{0}^{1}\alpha
                   \frac{d\beta}{dt}-\beta\frac{d\alpha}{dt}))\\
                  &=\text{exp}(\frac{1}{2\pi
                  i}\int_{0}^{1}(\log{m}\;d \log{l}-log{l}\;d
                  \log{m}))
   \end{aligned}
\end{equation}

where $(\alpha(t), \beta(t))$ is a lift of $\chi_t$ to
$\mathbb{C}^2$, under the $(l,m)$ coordinates,
$\displaystyle{\alpha=\frac{1}{2\pi i}\log{m}} $ and
$\displaystyle{\beta=\frac{1}{2\pi i}\log{l}}$ for a fixed branch
of logarithm.\\

It is clear that (\ref{kkf}) and (\ref{cs}) are not the same.
Moreover, the Chern-Simons $3$-form in \cite{KK} is the secondary
class from the second Chern class (a closed $4$-form). The term in
(\ref{cs}) defined in \cite[~(5.6)]{Guk} is a $1$-form. It may be
the secondary class of the first Chern class (a closed $2$-form)
of some line bundle over $\overline{D_K}$.\\

In the following subsections, we show that $dCS(l,m)$ indeed arises
from the first Chern class of a (universal) line bundle over the
Heisenberg group. Furthermore, we relate both $dVol$ and $dCS$ to
the imaginary and real parts of the secondary Chern-Simons class
respectively. We also give a mathematical proof of Gukov's
quantization statement of the Bohr-Sommerfield condition by using
some torsion element of $K_2$ and the regulator map.

\subsection{The regulator map of $K_2$}
Let $X$ be a smooth projective curve over $\mathbb{C}$ or a compact
Riemann surface. Let $f$, $g$ be two meromorphic functions on $X$.
Denote by $S(f)$ (resp. $S(g)$) the set of zeros and poles of $f$
(resp. $g$). Notice that $S(f)\cup S(g)$ is a finite set. Put
$X^{\prime}=X\setminus (S(f)\cup S(g))$.

Following Beilinson \cite{Bei}, see also \cite{De}, we define an
element $r(f,g)\in H^{1}(X^{\prime};\mathbb{C}^{*})$, equivalently,
as an element of $\text{Hom}(\pi_{1}(X^{\prime}),\mathbb{C}^{*})$:
for a loop $\gamma$ in $X^{\prime}$ with a distinguished base point
$t_0\in X^{\prime}$,
\begin{equation}\label{eq2.1}
r(f,g)(\gamma)=\exp{(\frac{1}{2\pi
i}(\int_{\gamma}\log{f}\;\frac{dg}{g}-\log{g(t_0)}\int_{\gamma}\frac{df}{f}))},
\end{equation}

where the integrals are taken over $\gamma$ beginning at $t_0$.

It is well-known that this definition is independent of the choices
of the base point $t_0$ and the branches of $\log{f}$ and $\log{g}$.
From now on, we shall take $\log{z}:\mathbb{C}^{*}\rightarrow
\mathbb{C}$ with $0\leq \arg{z}< 2 \pi$. Then it is well-defined,
but discontinuous on the positive real line $[0,+\infty)$ and it is
holomorphic on the cut plane $\mathbb{C}\setminus [0,+\infty)$.

In \cite{De}, Deligne noticed that
$H^{1}(X^{\prime};\mathbb{C}^{*})$ is the group of isomorphism
classes of the line bundles over $X^{\prime}$ with flat connections.
Hence $r(f,g)$ corresponds to such a line bundle with a flat
connection.
\begin{proposition}\label{prop1}
(1) The curvature of the line bundle associated to the class
$r(f,g)$ is
$\frac{df}{f}\wedge \frac{dg}{g}$; \\
(2) $r(f_{1}f_{2},g)=r(f_1,g)\otimes r(f_2,g)$,
$r(f,g)=r(g,f)^{-1}$,
and the Steinberg relation $r(f,1-f)=1$ holds if $f\ne 0$, $f\ne 1$;\\
(3) For $x\in S(f)\cup S(g)$, let $\gamma_{x}$ be a small simple
loop around $x$ in $X^{\prime}$. Then $r(f,g)(\gamma_{x})$ is equal
to the tame symbol $T_{x}(f,g)$ of $f$ and $g$ at $x$.
\end{proposition}

For the proof, we refer to \cite{De}. For the explicit construction
of the line bundle $r(f,g)$, see \cite[Section~4]{Ram} and \cite{Bl}
where the proof of this proposition was also given. The key
construction is a universal Heisenberg line bundle with connection
on $\mathbb{C}^{*}\times \mathbb{C}^{*}$. To prove the Steinberg
relation, the ubiquitous dilogarithm shows up.\\

Next recall the tame symbol
\[
   T_{x}(f,g):=(-1)^{v_{x}(f)\cdot
            v_{x}(g)}\frac{f^{v_{x}(g)}}{g^{v_{x}(f)}}(x),
\]
where $v_{x}(f)$ (resp. $v_{x}(g)$) is the order of zero or pole
of $f$ (resp. $g$) at $x$.\\

Let $\mathbb{C}(X)$ be the field of meromorphic functions on $X$.
Denote by $\mathbb{C}(X)^{*}$ the set of non-zero meromorphic
functions on $X$. By Matsumoto Theorem \cite{Mil},
\[
  K_2(\mathbb{C}(X))=\frac{\mathbb{C}(X)^{*}\otimes
               \mathbb{C}(X)^{*}}{\langle f\otimes (1-f): f\ne 0,1 \rangle},
\]

where the tensor product is taken over $\mathbb{Z}$, and the
denominator means the subgroup generated by those elements. For $f$
and $g\in \mathbb{C}(X)^{*}$, we denote by $\{f,g\}$ the
corresponding element in $K_2(\mathbb{C}(X))$.

The part $(2)$ of Proposition \ref{prop1} implies that we have a
homomorphism
\begin{equation}\label{eq2.2}
r: K_2(\mathbb{C}(X))\longrightarrow \underset{S\subset
X(\mathbb{C}):\;\text{finite}}{\varinjlim}H^{1}(X\setminus
S;\mathbb{C}^{*})
\end{equation}

defined by $r(\{f,g\})=r(f,g)$.

\subsection{}
Let $Y$ be an irreducible component of $\overline{D_K}$, the zero
locus of the $A$-polynomial $A(l,m)$. Denote by $\widetilde{Y}$ a
smooth projective model of $Y$. Then their fields of rational
functions are isomorphic,
$\mathbb{C}(Y)\cong\mathbb{C}(\widetilde{Y})$. We have the
following.
\begin{proposition}\label{prop2.1}
The element $\{l,m\}\in K_2(\mathbb{C}(Y))$ is a torsion element.
\end{proposition}

\begin{proof}
By \cite[Proposition~2.2, 4.1]{CCGLS}, there is a finite field
extension $F$ of $\mathbb{C}(Y)$ such that $\{l,m\}\in K_2(F)$ is of
order at most $2$. We have a homomorphism $i:
K_2(\mathbb{C}(Y))\rightarrow K_2(F)$ induced by the inclusion of
$\mathbb{C}(Y)$ into $F$. We also have the transfer map
$t:K_2(F)\rightarrow K_2(\mathbb{C}(Y))$. It is well-known that the
composition $t\circ i$:
\[
  K_2(\mathbb{C}(Y))\rightarrow K_2(F)\rightarrow K_2(\mathbb{C}(Y))
\]
is given by the multiplication of $n=[F:\mathbb{C}(Y)]$, the degree
of the finite extension. Hence $t(i(\{l,m\})=t(\{l,m\})=n\{l,m\}$.
This implies that $\{l,m\}\in K_2(\mathbb{C}(Y))$ is a torsion and
its order divides $2n$.
\end{proof}

Suppose the component $Y$ contains $y_0\in D_K$ which corresponds to
the discrete faithful character $\chi_0$ of the hyperbolic structure
and $m(y_0)=1$. Let $S(l,m)$ be the finite set of poles and zeros of
$l$ and $m$. Put $Y_h=\widetilde{Y}\setminus S(l,m)$ as the
$X^{\prime}$ in \S 3.2. We choose the distinguished point $t_0$ as
follows. If $y_0$ is a smooth point, we take $t_0=y_0$; if $y_0$ is
a singular point, we fix a point in the pre-images of $y_0$ in
$\widetilde{Y}$ and take $t_0$ as this fixed point. This is
equivalent to fixing a branch at the singular point $y_0$.
\begin{theorem}\label{main}
(i) The closed real $1$-form $\eta(l,m)=\log{|l|}\;
d\arg{m}-\log{|m|}\; d\arg{l}$ is exact on $Y_h$;\\
(ii) For any loop $\gamma$ with initial point $t_0=\chi_0$ in $Y_h$
\[
   \frac{1}{4
   \pi^{2}}\int_{\gamma}(\log{|m|}\, d\log{|l|}+\arg{l}\, d\arg{m})=
    \frac{p}{q},
\]
where $p$ is some integer and $q$ is the order of the symbol
$\{l,m\}$ in $K_2(\mathbb{C}(Y))$.
\end{theorem}

\begin{proof}
By (\ref{eq2.2}), we have an element $r(l,m)\in
H^{1}(Y_h,\mathbb{C}^{*})$. By Proposition \ref{prop2.1}, it is a
torsion of order $q$. By the definition of $r(l,m)$ in
(\ref{eq2.1}), we conclude that for any loop $\gamma$ in $Y_h$,
\begin{equation}\label{kt}
    \{\exp{(\frac{1}{2\pi
    i}(\int_{\gamma}\log{l}\;\frac{dm}{m}-\log{m(t_0)}\int_{\gamma}\frac{dl}{l}))}\}^q=1
\end{equation}

Write
$\displaystyle{\int_{\gamma}\log{l}\;\frac{dm}{m}-\log{m(t_0)}\int_{\gamma}\frac{dl}{l}}=Re+iIm$,
where $Re$ and $Im$ are the real and imaginary parts respectively.
(\ref{kt}) means that $\displaystyle{\exp{(\frac{q\cdot
Im}{2\pi}+\frac{q\cdot Re}{2\pi i})}=1}$. Therefore, $Im=0$ and
$\displaystyle{\frac{q\cdot Re}{2\pi i}=2 \pi i p}$, for some
integer $p$. Our result follows from the following lemma.
\end{proof}

\begin{lemma}\label{le2.1}
Denote
$\int_{\gamma}\log{l}\;\frac{dm}{m}-\log{m(t_0)}\int_{\gamma}\frac{dl}{l}=Re+iIm$
as above. Then
\[
Im=\int_{\gamma}(\log{|l|}\,d\arg{m}-\log{|m|}\,d\arg{l})=\int_{\gamma}\eta(l,m),
\]
and
\[
Re=-\int_{\gamma}(\log{|m|}\,d\log{|l|}+\arg{l}\,d\arg{m})=\int_{\gamma}\xi(l,m),
\]
where $\xi(l,m)$ depends on the branches of $\arg$ function and $Re$
is well-defined up to $(2\pi)^{2}\mathbb{Z}$.
\end{lemma}

\begin{proof}
Let $F$ be a smooth non-zero complex-valued function, and
$F=Re(F)+iIm(F)$, where $Re(F)$ denotes its real part and $Im(F)$
its imaginary part. Then we have
\[
d\log{F}:=\frac{dF}{F}=\frac{d|F|}{|F|}+i\frac{Re(F)dIm(F)-Im(F)dRe(F)}{|F|^2}
\]

So the real part of $d\log{F}$ is $d\log{|F|}$ which is exact and
the imaginary part is denoted by $d\arg{F}$.

By a straightforward calculation, we have
\[
Im=\int_{\gamma}(\log{|l|}d\arg{m}+\arg{l}d\log{|m|})-\log{|m(t_0)|}\int_{\gamma}d\arg{l}.
\]

Integration by parts, we obtain:
\[
\int_{\gamma}\arg{l}\cdot
d\log{|m|}=\log{|m(t_0)|}\int_{\gamma}d\arg{l}-\int_{\gamma}\log{|m|}\cdot
d\arg{l}.
\]
Therefore,
\[
  Im=\int_{\gamma}(\log{|l|}d\arg{m}-\log{|m|}d\arg{l}).
\]

   For the real part $Re$, it is equal to
\[
\int_{\gamma}(\log|l|\;d\log|m|-\arg{l}\;d\arg{m})+\arg{m(t_0)}\int_{\gamma}d\arg{l}.
\]

Integration by parts, we get
\[
\int_{\gamma}\log|l|\,d\log|m|=-\int_{\gamma}\log{|m|}\,d\log{|l|}.
\]
By the choice of $t_0=\chi_{\rho_0}$, $\arg{m(t_0)}=0$. Hence the
result follows.
\end{proof}

\begin{remark}
(i) The first part of the theorem was also proved in
\cite[Sect.~4.2]{CCGLS}. Our proof gives the real and imaginary part
simultaneously.

(ii) The result of the second part is stronger than the one in \cite
[3.29]{Guk} where he derived that the value of the integral is in
$\mathbb{Q}$ from the quantizable Bohr-Sommerfield condition.

(iii) The class $r(l,m)\in H^1(Y_h; \mathbb{C}^{*})$ corresponds to
a flat line bundle $L$ over $Y_h$ which is the pullback of the
universal Heisenberg line bundle on $\mathbb{C}^{*}\times
\mathbb{C}^{*}$, see \cite{Bl,Ram}. Formally,
\[
  d(\xi(l,m)+i\eta(l,m))=\frac{dl}{l}\wedge \frac{dm}{m}=0.
\]
Hence, $\frac{1}{2 \pi i}(\xi(l,m)+i\eta(l,m))$ is the $1$-form
Chern-Simons. Denote it by $CS_1(l,m)$. Then
$dCS_1(l,m)=C_1(L)=\frac{1}{2 \pi i}\frac{dl}{l}\wedge
\frac{dm}{m}=0$.
\end{remark}

By Theorem \ref{main} and Lemma \ref{le2.1}, we would like to
propose the corresponding generalized volume conjecture as the
following:\\

For a path $c: [0,1]\rightarrow Y_h$ with $c(0)=t_0$ and
$c(1)=(l,m)$, write $c(t)=(l(t),m(t))$. Recall that $q$ is the order
of the symbol $\{l,m\}$ in $K_2(\mathbb{C}(Y))$. We denote
\[
  U(l,m)=-q \cdot \int_{c}[\log{|m(t)|}\,d\log{|l(t)|}+\arg{l(t)}\,d\arg{m(t)}].
\]
and call it the \emph{special Chern-Simons invariant} of $K$.

For a fixed number $a$ and $m=-\exp({i \pi a})$, we reformulate the
generalized volume conjecture as the following:\\

\noindent{\bf Conjecture}: (The reformulated generalized volume
conjecture)
\begin{equation}\label{re-conj}
\lim_{N, k\rightarrow \infty; \frac{N}{k}=a}
\frac{\log{J_{N}(K,e^{2\pi i/k})}}{k}=\frac{1}{2 \pi}(Vol(l,m)+
i\frac{1}{2\pi}U(l,m)).\\
\end{equation}
\newline

\begin{remark}
By Theorem \ref{main} (ii), $\frac{1}{(2\pi)^2}U(l,m)$ is
well-defined in $\mathbb{R}/\mathbb{Z}$. The classical Chern-Simons
invariant is also well-defined in $\mathbb{R}/\mathbb{Z}$.
\end{remark}

Our reformulated generalization (\ref{re-conj}) gives a
$\mathbb{C}^{*}$-paramatrized version of the volume conjecture.
Using Fuglede-Kadison determinant, W. Zhang and the first author
\cite{LZ} defined an $L^2$-version twisted Alexander polynomial
which can be identified with an $L^2$-Reidemeister torsion. By Luck
and Schick's result, this $L^2$-Alexander polynomial provides a
$\mathbb{C}^{*}$-parametrization of the hyperbolic volumes. For
other discussions on the volume conjecture, see
\cite[Sections~1.3,\;7.3]{Oh} and the related references within.

\begin{remark}
In \cite{GuMu}, Gukov and Murakami showed that the difference of
their conjectures (in \cite{Guk} and \cite{Mu} respectively) comes
from the different choice of polarization. One choice leads to
Gukov's, the other gives Murakami's. Our generalization
(\ref{re-conj}) only focuses on Gukov's original one in \cite{Guk}.
Moveover, it indicates that the special Chern-Simons term $U(l,m)$
comes from the regulator over the character variety. It would be
interesting to extend our reformulated generalized conjecture
(\ref{re-conj}) to other polarizations.
\end{remark}

\subsection{A Motivic Perspective}
In \cite{Gon}, A. Goncharov proved that the volume of a hyperbolic
$3$-manifold is the period of a mixed Tate motive, and he gave the
explicit construction of this motive there. The Hopf algebra
$\mathcal {H}$ of framed mixed Tate motives was also defined there,
for the detail, see \cite{Gon} and the references therein.

For our case, the equation (\ref{vol}) can be thought of as the
variation formula of the volume on the deformation space $X_0$ of
hyperbolic structures on $M_K$, see \cite{Th1, Hod, CCGLS}. On the
motivic level, when we deform the hyperbolic structures on $M_K$, we
can think of it as the variation of the corresponding mixed motives.
Professor A. Goncharov pointed out to the second-name author that in
this picture, the equation (\ref{vol}) is related to the coproduct
of the Hopf algebra $\mathcal {H}$. This indicates an interesting
link for a future study.


\begin{thebibliography}{99}
\bibitem[Bei]{Bei} Beilinson, A., \textit{Higher regulators and
values of $L$-functions of curves}, Funktsional. Anal. i Prilozhen.
14 (1980), no. 2, 46--47.
\bibitem[Bl]{Bl} Bloch, S., \textit{The dilogarithm and extensions
of Lie algebras}, Algebraic $K$-theory, Evanston 1980, Lecture Notes
in Math., 854, 1-23, Springer, Berlin-New York, 1981.
\bibitem [CCGLS]{CCGLS} Cooper, D., Culler, M., Gillet, H.,
Long, D.D., Shalen, P.B., \textit{Plane curves associated to
character varieties of 3-manifolds}, Invent. Math. 118(1994), 47-74.
\bibitem [CGLS]{CGLS} Culler, M., Gordon, C., Lueke, J., Shalen, P.,
\textit{Dehn Surgery on Knots}, Ann. of Math. 125(1987), 297-930.
\bibitem [CS]{CS} Culler, M., Shalen, P.B., \textit{Varieties of
group representations and splittings of 3-manifolds}, Ann. Math.,
(2) 117(1464), 109-146.
\bibitem [De]{De} Deligne, P., \textit{Le symbole
mod\'{e}r\'{e}}, Inst. Hautes ¨¦tudes Sci. Publ. Math. No. 73,
(1991), 147--181.
\bibitem [Dun]{Dun} Dunfield, M.N., \textit{Cyclic surgery, degrees
of maps of character curves, and volume rigidity for hyperbolic
manifolds}, Invent. Math. 170(1999) 623-650.
\bibitem [Gon]{Gon} Goncharov, A., \textit{Volumes of hyperbolic manifolds and mixed Tate
motives},  J. Amer. Math. Soc.  12  (1999),  no. 2, 569--618.
\bibitem [Guk]{Guk} Gukov, S., \textit{Three-dimensional quantum
gravity, Chern-Simons theory, and the A-polynomial},
Comvun.Math.Phys. 055(2005), 577-622.
\bibitem [GuMu]{GuMu} Gukov, S., Murakami, H., \textit{SL(2,C) Chern-Simons theory and the asymptotic behavior of the colored Jones
polynomial}, arXiv:math.GT/0608324.
\bibitem [Hod]{Hod} Hodgson, C., \textit{Degeneration and regeneration of hyperbolic structures on three-manifolds}, Thesis,
Princeton University, 1986.
\bibitem [Kas]{Kas} Kashaev, R., \textit{The hyperbolic volume of
knots from quantum dilogarithm}, Lett. Math. Phys. 39(1997),
269-275.
\bibitem [KK]{KK} Kirk, P., Klassen, E., \textit{Chern-Simons
Invariants of 3-manifolds decomposed along tori and the circle
bundle over the representation space of $T^2$}, Commun. Math. Phys.
153, 581-557(1393).
\bibitem[LZ]{LZ} Li, W., Zhang, W., \textit{An $L^2$-Alexander
invariant for knots}, Commun. Contemporary Math., Vol.8, No.2,
(2006) 1-21.
\bibitem[Mil]{Mil} Milnor, J., \textit{Introduction to algebraic $K$-theory},
Annals Math. Studies, No. 72. Princeton University Press, 1971.
\bibitem[MM]{MM} Murakami, H., Murakami, J., \textit{The colored
Jones polynomial and the simplicial volume of a knot}, Acta. Math.
186(2001), 85-104.
\bibitem[Mu]{Mu} Murakami, H., \textit{A version of the volume
conjecture}, arXiv:math.GT/0603217.
\bibitem[Oh]{Oh} Ohtsuki, T., \textit{Problems on invariants of
knots and $3$-manifolds}, Goemetry and Topology Monographs, Vol 4,
377-572.
\bibitem [Ram]{Ram} Ramakrishnan, D., \textit{Regulators, algebraic
cycles, and values of L-functions}, Contemporary Math. 83,
183-310(1989).
\bibitem [RSW]{RSW} Ramadas, T., Singer, I., Weitsman, J.,
\textit{Some comments on Chern-Simons gauge theory},
Commun.Math.Phys. 126(1989), 409-420.
\bibitem [Sha]{Sha} Shalen, P.B., \textit{Representations of hyperbolic manifolds}, In: \textit{Handbook of geometric topology},
Elsevier Press.
\bibitem [Th1]{Th1} Thurston, W., \textit{The geometry and topology of 3-manifolds}, Lecture Notes, Princeton Univ. Press, 1978.
\bibitem [Th2]{Th2} Thurston, W., \textit{Three-dimensional
manifolds, Kleinian groups and Hyperbolic geometry}, Bull. Amer.
Math. Soc. (N.S.) 6, 357-381(1982).


\end{thebibliography}
\end{document}